\documentclass[]{article}
\usepackage{graphicx}
\usepackage{caption}
\usepackage{algorithmic}
\usepackage{algorithm}
\usepackage{amsmath}
\usepackage{amssymb}
\usepackage{amsfonts}
\usepackage{amsthm}
\usepackage{fullpage}
\usepackage{url}
\usepackage{color}
\usepackage{authblk} 
\usepackage{float}
\usepackage{mathabx} 
\usepackage{subfig}

\def\N{{\mathbb{N}}}
\def\R{{\mathbb{R}}}

\def\err{{\textrm{err}}}
\def\IMF{{\textrm{IMF}}}
\def\dFIF{{\textrm{dFIF}}}
\def\htFIF{{\textrm{htFIF}}}
\def\FIF{{\textrm{FIF}}}

\newtheorem{theorem}{Theorem}

\title{Iterative Filtering as a direct method for the decomposition of non--stationary signals}
\author{Antonio Cicone}

\begin{document}

\maketitle

\begin{abstract}
The Iterative Filtering method is a technique developed recently for the decomposition and analysis of non-stationary and non-linear signals.
In this work we propose two alternative formulations of the original algorithm which allows to transform the Iterative Filtering method into a direct technique, making the algorithm closer to an online algorithm. We present a few numerical examples to show the effectiveness of the proposed approaches.
\end{abstract}

\section{Introduction}\label{sec:Intro}

The Iterative Filtering (IF) method is, as the name suggests, an iterative algorithm proposed by Lin et al. in 2009 \cite{lin2009iterative} as an alternative to the well known Empirical Mode Decomposition (EMD) method. The EMD is part of the so called Hilbert Huang Transform (HHT) technique \cite{huang1998empirical} for the analysis of non-stationary and non-linear signals. The aim of the EMD and IF is the decomposition of a given signal into simple components, defined by Huang \cite{huang1998empirical} as Intrinsic Mode Functions (IMFs), which are oscillatory functions that fulfills two empirical properties: the two curves connecting respectively the maxima and minima of an IMF have to be symmetric with respect to the horizontal line; the number of zero crossing must equal plus or minus one the number of local extrema of an IMF \cite{cicone2017dummies}.

Once a non-stationary and non-linear signal has been decomposed into IMFs with different scales it is possible to study its properties to unravel potentially hidden features.  Furthermore such decomposition allows for a more accurate time frequency analysis of the signal itself. The EMD and IF decomposition methods have been applied to the study of a wide variety of datasets. For instance in Medicine they have been used for the automatic identification of seizure-free electroencephalographic signals \cite{pachori2008discrimination}, for the study of the gastroesophageal reflux disease  \cite{liang2005application}, for the derivation of the respiratory sinus arrhythmia from the heartbeat time series \cite{balocchi2004deriving}, for the analysis of the heart rate variability \cite{echeverria2001application}, to denoise the electrocardiographic signal and correct the baseline wander \cite{blanco2008ecg}, to study the dengue virus spread \cite{gubler2004cities}, for the identification of the coupling between prefrontal and visual cortex \cite{gregoriou2009high}; in Geophysics, for the study of the evolution of land surface air temperature trend \cite{ji2014evolution}, to analyze the global mean sea-level rise \cite{chen2017increasing}, for the extraction of the solar cycle signal from the stratospheric data \cite{coughlin200411}, to identify near-fault ground-motion characteristics \cite{loh2001application}, for the study of electromagnetic field observations from satellites during earthquakes \cite{piersanti2018AQ}; in Engineering, to diagnose faults in rotating machinery \cite{lei2013review}, to separate two sources from a single-channel recording \cite{mijovic2010source}, to control the wind response of a building with variable stiffness tuned mass damper \cite{varadarajan2004wind}, to improve the detection of thermal bridges in buildings \cite{sfarra2019thermal}, to improve the detection of chemical plumes \cite{cicone2016hyperspectral}; in Information Technology, to analyze images \cite{nunes2003image} and texture \cite{nunes2005texture}; in Economics, to analyze the price of the crude oil \cite{zhang2008new}. This seemingly long list of applications is actually far from being complete. Consider, in fact, that the first paper ever published on the EMD method \cite{huang1998empirical} has received so far by itself, based on the Scopus database, more than 10800 citations.

In a recent work \cite{cicone2019IF_num_an} the complete numerical analysis of the IF has been addressed. One of the consequences of this analysis is that the IF algorithm convergence is guaranteed a priori and, assuming a periodical boundary extension of a signal, the IF can be implemented on a computer using the Fast Fourier Transform (FFT). The derived method, called Fast Iterative Filtering (FIF), is consistently faster than EMD and it allows to reformulate the IF method to become a direct algorithm. We recall that the study of the IF method for more general boundary conditions as been tackled in \cite{cicone2017BC}. Furthermore the extension of IF to higher dimensions has been studied in \cite{cicone2017multidimensional}. Finally it is important to remind here that the IF algorithm can be generalized to become the so called Adaptive Local iterative Filtering (ALIF) method which allows to better identify chirps contained in a signal. However the analysis of the ALIF technique is still in progress \cite{cicone2016adaptive,cicone2017Geophysics,cicone2017spectral}.

The rest of the paper is organized as follows. In Section \ref{sec:Hard_soft} we propose two new alternative formulations of IF as a direct technique. In Section \ref{sec:examples} we study a few numerical examples to show the abilities of these newly proposed methods compared with the original method. We conclude this work with a few comments.

\section{Iterative Filtering as a direct method}\label{sec:Hard_soft}

We start this section by recalling the basic Iterative Filtering (IF) method when applied to a discrete signal $s(x)$, $x\in\R$. We assume for simplicity that $s$ is supported on $[0,\ 1]$, sampled at $n$ points $x_j= \frac{j}{n-1}$, with $j= 0,\ldots, n-1$, and a sampling rate which allows to capture all its fine details, so that aliasing will not play any role. The goal of the  algorithm is to decompose $s$ into IMFs. Given the matrix
\begin{equation}\label{eq:K}
    W=\left[w(x_i-x_j)\cdot \frac{1}{n}\right]_{i,\ j=0}^{n-1}
\end{equation}
the main step of IF is given by
\begin{equation}\label{eq:IF_MainStep}
    s_{m+1} = (I-W)s_m \qquad m\in\N
\end{equation}
where $s_1=s$, and $w(x_j)$ is a nonnegative, even, and compactly supported window/filter with area equal to one and support in $[-L,\ L]$. $L$ is called filter length and is calculated based on the signal $s$ itself \cite{cicone2019IF_num_an}.

We recall also the following
\begin{theorem}[Convergence of IF \cite{cicone2019IF_num_an}]
    Given a discrete signal $s$, a filter function $w$ whose support length is based on the signal itself and assuming to extend periodically at the boundaries the signal, then the first IMF is given by
\begin{equation}\label{eq:IMF1_direct}
\IMF_1=U(I-D)^{N_0} U^T s = \textrm{IDFT}\left((I-\textrm{diag}\left(\textrm{DFT}(w)\right))^{N_0}\textrm{DFT}(s)\right)
\end{equation}
where $D$ is a diagonal matrix containing on the diagonal the eigenvalues of $W$, $U$ is matrix having as columns the eigenvectors of the circulant matrix $W$, and $N_0$ is the number of iterations needed to compute the first IMF, based on a predefined stopping criterion.
\end{theorem}
We point out that subsequent IMFs can be computed using the same formula applied to the remainder $s-\sum_{k=1}^M\IMF_k$, where $M$ represents the number of IMFs already computed.

Thanks to this theorem the computation of $s_{m+1}$ can be made fast using the Fast Fourier Transform (FFT). This is the idea behind what is called the Fast Iterative Filtering (FIF) algorithm \cite{cicone2019IF_num_an}.

From \eqref{eq:IMF1_direct} it is clear that the role of the iterations is that of sending to zero all the eigenvalues of $I-D$ which are not equal to 1 and preserving unchanged the ones equal to 1.
Therefore  we can make the iterative method to become a direct one by properly transforming the eigenvalues of $I-D$.
This is the idea behind both the direct Fast Iterative Filtering (dFIF) and the hard thresholding Fast Iterative Filtering (htFIF) that we are going to present.

\subsection{Direct Fast Iterative Filtering}

The \emph{Direct Fast Iterative Filtering} (dFIF) works as follows: given the eigenvalues $\lambda_i\in\sigma(W)$, for $i=1,\ \ldots,\ n$, we set a threshold $\tau$ and a value $\kappa=0.5$, so that
we can estimate an appropriate value for $N_0$ which allows to compute each IMF in one step. In particular we choose
\begin{equation}\label{eq:N0}
    N_0 = \textrm{round} \left(\frac{\log\left(\kappa\right)}{\log\left(\displaystyle \max_{1-\lambda_i<\tau}\left(1-\lambda_i\right)\right)}\right).
\end{equation}

The pseudo code of the dFIF method is given in Algorithm \ref{algo:dFIF}

\begin{algorithm}
\caption{\textbf{Direct Fast Iterative Filtering} IMF = dFIF$(s)$}\label{algo:dFIF}
\begin{algorithmic}
\STATE IMF = $\left\{\right\}$
\WHILE{the number of extrema of $s$ $\geq 2$}
      \STATE compute the filter function $w$
      \STATE compute $N_0$ using eqref{eq:N0}
      \STATE $\IMF_k = \textrm{IDFT}\left((I-\textrm{diag}\left(\textrm{DFT}(w)\right))^{N_0}\textrm{DFT}(s)\right)$
      \STATE IMF = IMF$\,\cup\,  \{ \IMF_{k}\}$
      \STATE $s=s-\IMF_k$
\ENDWHILE
\STATE IMF = IMF$\,\cup\,  \{ s\}$
\end{algorithmic}
\end{algorithm}

The numerical relative error introduced by this approach to compute the first IMF is given by
\begin{equation}\label{eq:error_dFIF}
\textrm{err}_{\dFIF} = \frac{\left\|\IMF_{\dFIF}-\IMF_{\FIF}\right\|_2}{\left\|\IMF_{\FIF}\right\|_2} \leq  \frac{\left\|\left((I-\textrm{diag}\left(U w\right))^{N_{\dFIF}-N_{\FIF}}-I\right)(I-\textrm{diag}\left(U w\right))^{N_{\FIF}}\right\|_2\left\|s\right\|_2}{\left\|\left(I-\textrm{diag}\left(U w\right)\right)^{N_{\FIF}}U^T s\right\|_2}
\end{equation}
where we assume, for simplicity, that $N_{\dFIF}>N_{\FIF}$.

The error goes to zero if the value $N_{\dFIF}$, estimated using \eqref{eq:N0}, equals the number of steps $N_{\FIF}$ needed by the FIF algorithm to compute the first IMF.

We point out that in \cite{cicone2019IF_num_an} the authors proposed another formula for the a priori computation of $N_0$. However such formula provides only a rough overestimate of the actual number $N_0$, whereas \eqref{eq:N0} can give a better estimate for the value $N_0$.

The dFIF has the advantage of speeding up the calculations with respect to the FIF method and potentially having a zero $\textrm{err}_{\dFIF}$ error.

However it is important to point out here that the error we are computing is not with respect to the exact solution. The FIF decomposition represents by itself only an approximation of the exact one. Therefore having a zero $\textrm{err}_{\dFIF}$ is not necessarily the perfect measure of the quality of the decomposition.

\subsection{Hard Thresholding}

The other possible approach is called \emph{Hard Thresholding Fast Iterative Filtering} (htFIF): we set to zero all the eigenvalues of the matrix $I-D$ which are smaller than a threshold $\tau$, and we leave unchanged all the other ones.

The pseudo code of the htFIF is given in Algorithm \ref{algo:htFIF}

\begin{algorithm}
\caption{\textbf{Hard Thresholding Fast Iterative Filtering} IMF = htFIF$(s)$}\label{algo:htFIF}
\begin{algorithmic}
\STATE IMF = $\left\{\right\}$
\WHILE{the number of extrema of $s$ $\geq 2$}
      \STATE compute the filter function $w$
      \STATE $A = I-\textrm{diag}\left(\textrm{DFT}(w)\right)$
      \STATE $B$ is equal to $A$ with all the diagonal entries smaller than $\tau$ set to zero
      \STATE $\IMF_k = \textrm{IDFT}\left(B^{N_0}\textrm{DFT}(s)\right)$
      \STATE IMF = IMF$\,\cup\,  \{ \IMF_{k}\}$
      \STATE $s=s-\IMF_k$
\ENDWHILE
\STATE IMF = IMF$\,\cup\,  \{ s\}$
\end{algorithmic}
\end{algorithm}

The numerical relative error introduced by this second approach is given by
\begin{equation}\label{eq:error_htFIF}
 \textrm{err}_{\htFIF} = \frac{\left\|\IMF_{\htFIF}-\IMF_{\FIF}\right\|_2}{\left\|\IMF_{\FIF}\right\|_2} \leq \frac{\left\|B-(I-\textrm{diag}\left(U w\right))^{N_{\FIF}}\right\|_2 \left\|s\right\|_2}{\left\|\left(I-\textrm{diag}\left(U w\right)\right)^{N_{\FIF}}U^T s\right\|_2}
\end{equation}

The advantage of htFIF with respect to the dFIF is that we have only one parameter to tune. In this case the value of $\textrm{err}_{\htFIF}$ can not become zero in general and the decomposition we expect to be, potentially, less accurate than the ones produced using either dFIF or FIF methods.

\section{Numerical examples}\label{sec:examples}

In this section we run comparisons of the decompositions produced via FIF, dFIF and htFIF.
We do not make any comparison with other iterative methods available in the literature, i.e. the Empirical Mode Decomposition method and the Ensemble Empirical Mode Decomposition because the FIF method is already by itself faster than both of them.

All the calculations presented in this section have been performed on a laptop equipped with an Intel Core i7-8550U CPU, 1.80GHz, 16.0 GB RAM, Windows 10 Pro, Matlab R2018a. The FIF code is available at \url{http://www.cicone.com}.

\subsection{Example 1}
\begin{figure}%
    \centering
    \subfloat{{\includegraphics[width=0.48\textwidth]{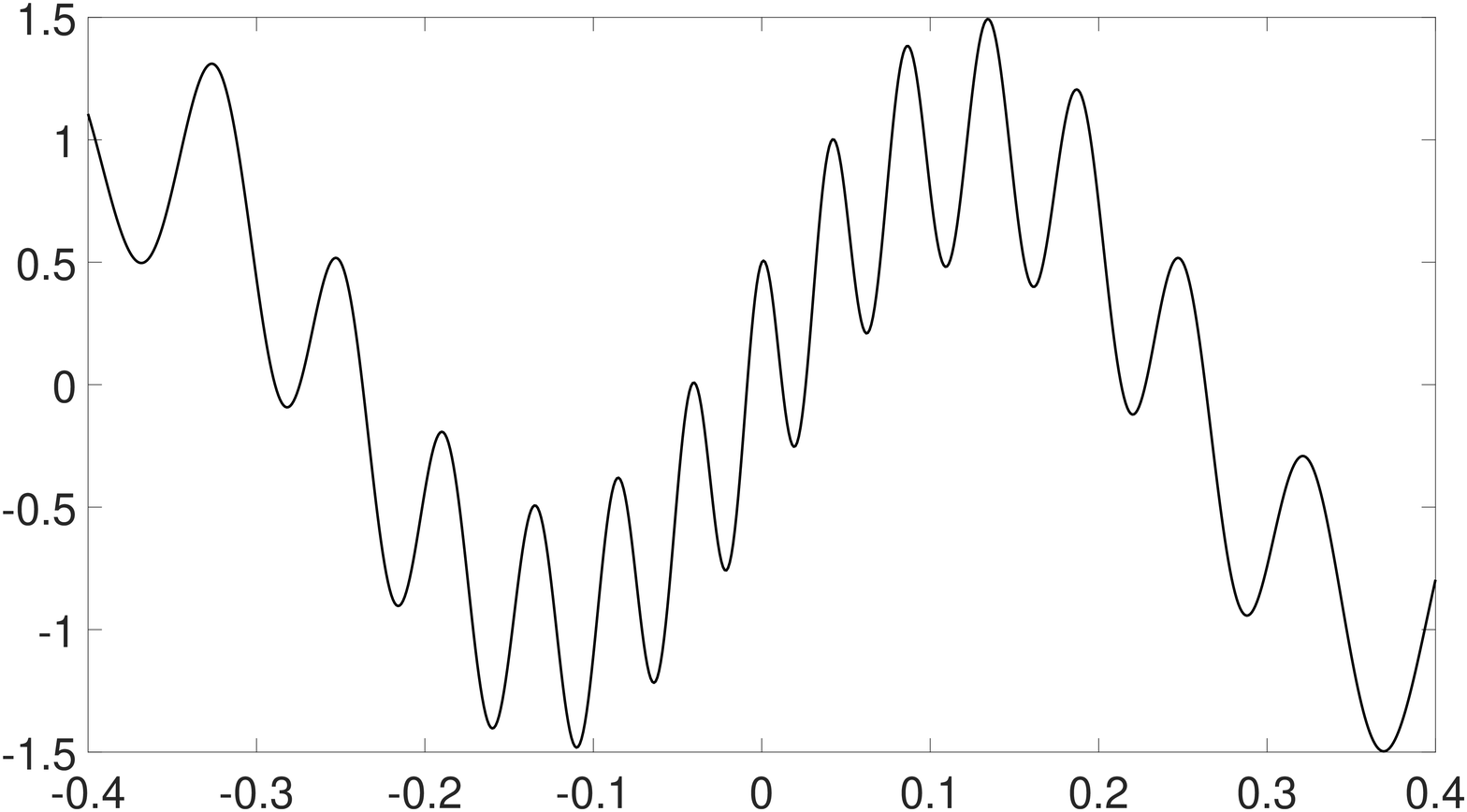} }}%
    ~
    \subfloat{{\includegraphics[width=0.48\textwidth]{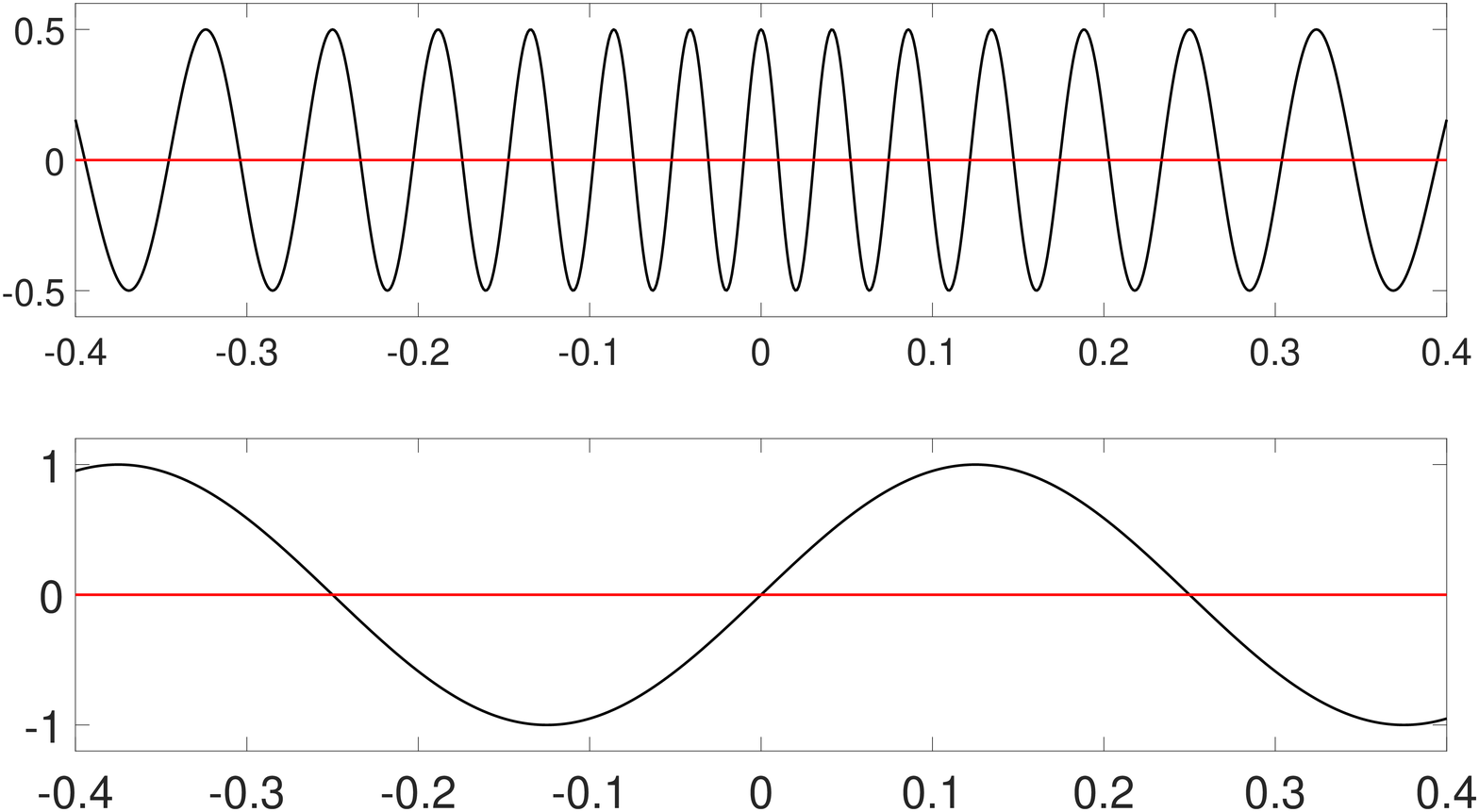} }}%
    \caption{ Left panel, artificial signal. Right panel, ground truth}\label{fig:Test_ex_2}
\end{figure}

We start with the toy example shown in Figure \ref{fig:Test_ex_2} which is produced as the summation of the two oscillatory signals plotted the same figure.

\begin{figure}%
    \centering
    \subfloat{{\includegraphics[width=0.32\textwidth]{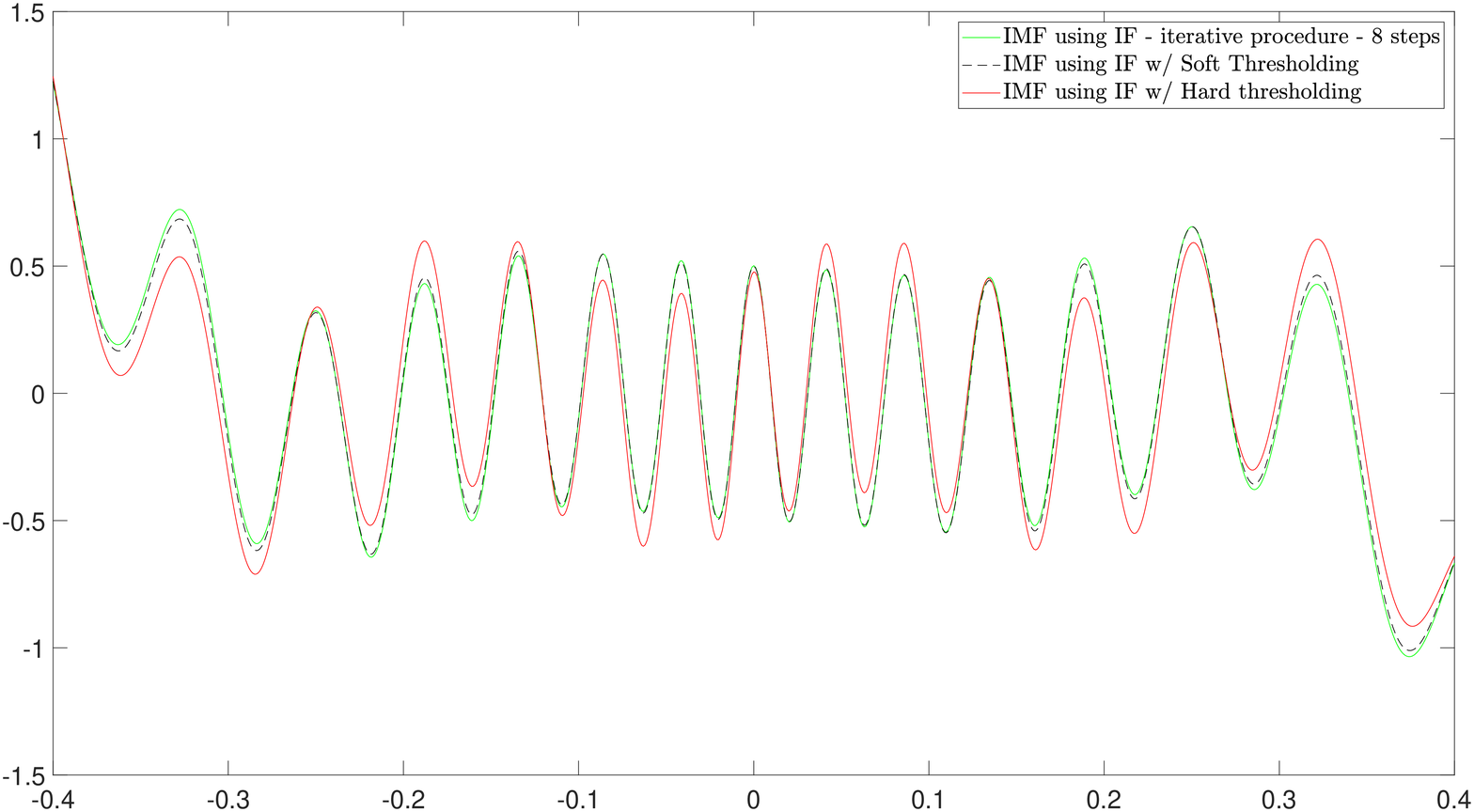} }}%
    ~
    \subfloat{{\includegraphics[width=0.32\textwidth]{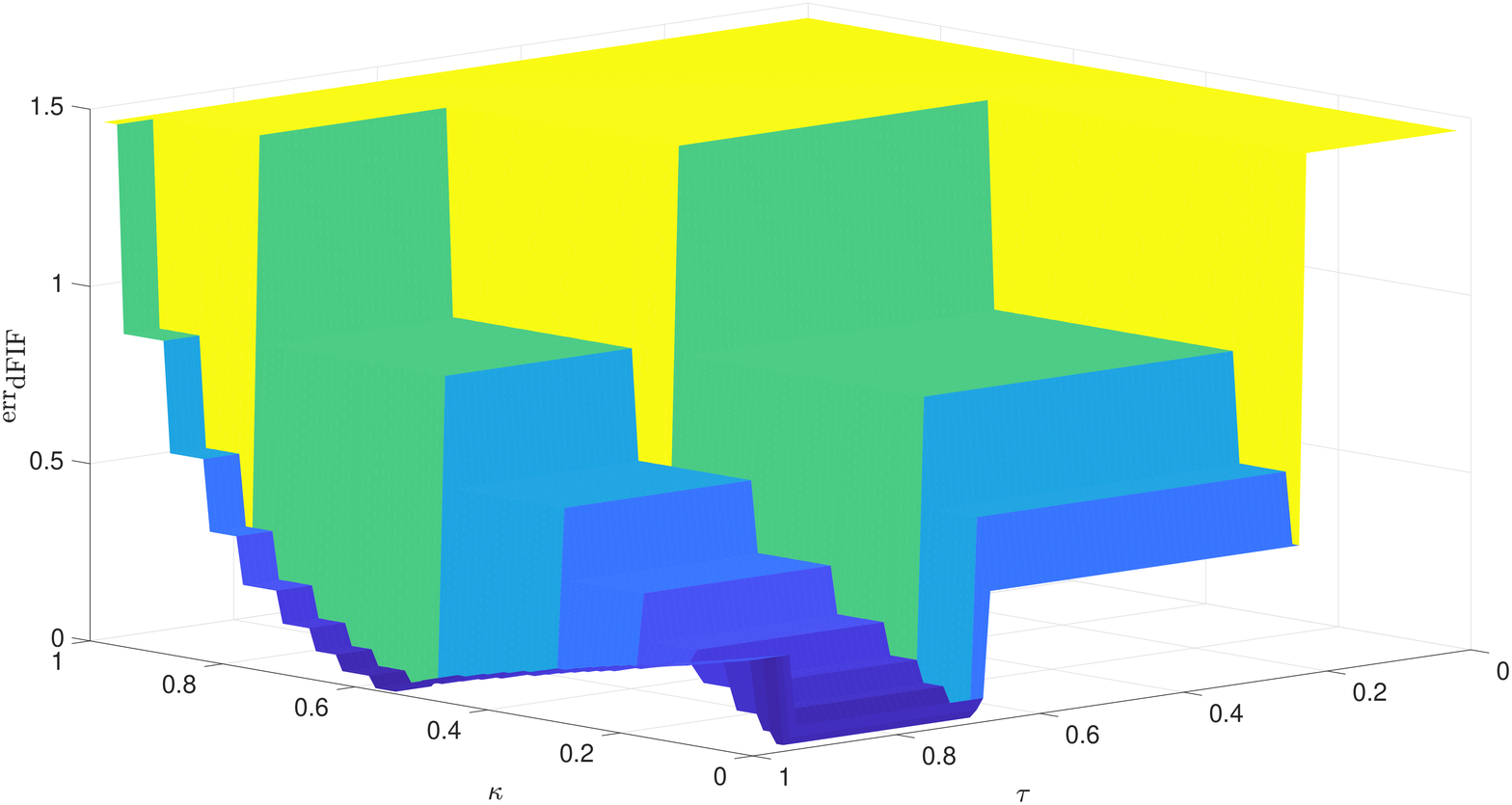} }}%
    ~
    \subfloat{{\includegraphics[width=0.32\textwidth]{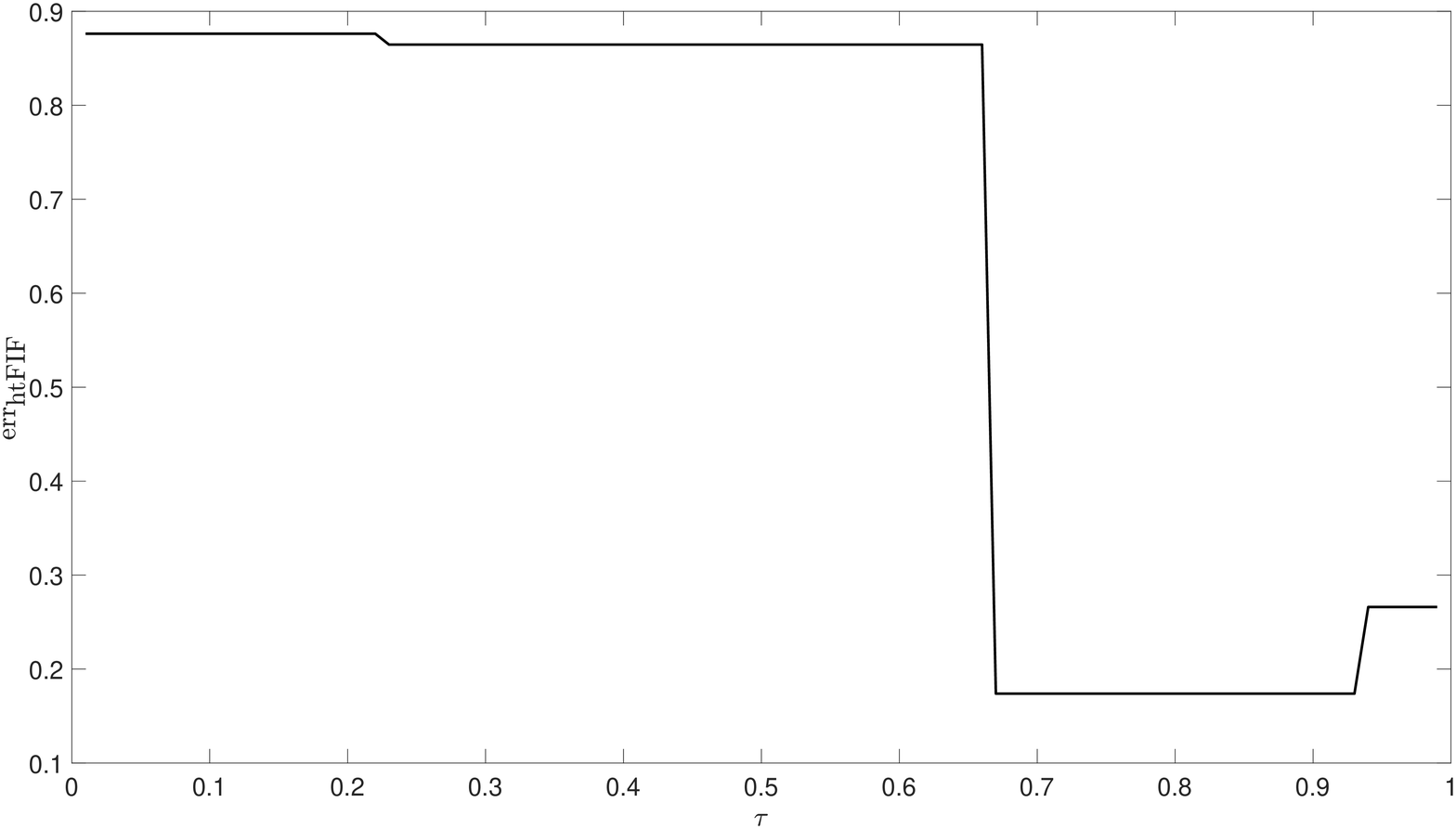} }}%
    \caption{ Left panel, first IMF produced using the standard FIF algorithm, FIF based on Soft and Hard thresholding. Middle panel, $\err_\dFIF$ values for different values of $\tau$ and $\kappa$. Right panel, $\err_\htFIF$ values for different values of $\tau$.}\label{fig:Test_ex_2_errors}
\end{figure}

In the left panel of Figure \ref{fig:Test_ex_2_errors} we compare the first IMFs produced by the FIF\footnote{We set \textbf{alpha} equal to 1, \textbf{Xi} $=4$, \textbf{delta} $=0.001$.}, dFIF, and htFIF methods. We compute the relative errors of dFIF and htFIF for different values of $\tau$ and $\kappa$. The results are shown in the middle and right panel of Figure \ref{fig:Test_ex_2_errors}. From this analysis we see that the optimal $\tau$ for htFIF is in the interval $[0.7,\ 0.9]$. Whereas from the dFIF method results it follows that either we select $\tau>0.95$ and $\kappa\approx 0.56$ or $\tau\in[0.65,\ 0.95]$ and $\kappa\approx 0.02$.

Based on these results we decide to set $\tau=0.98$ and $\kappa= 0.56$ in the subsequent examples.

\subsection{Example 2}

\begin{figure}
    \centering
    \subfloat{{\includegraphics[width=0.48\textwidth]{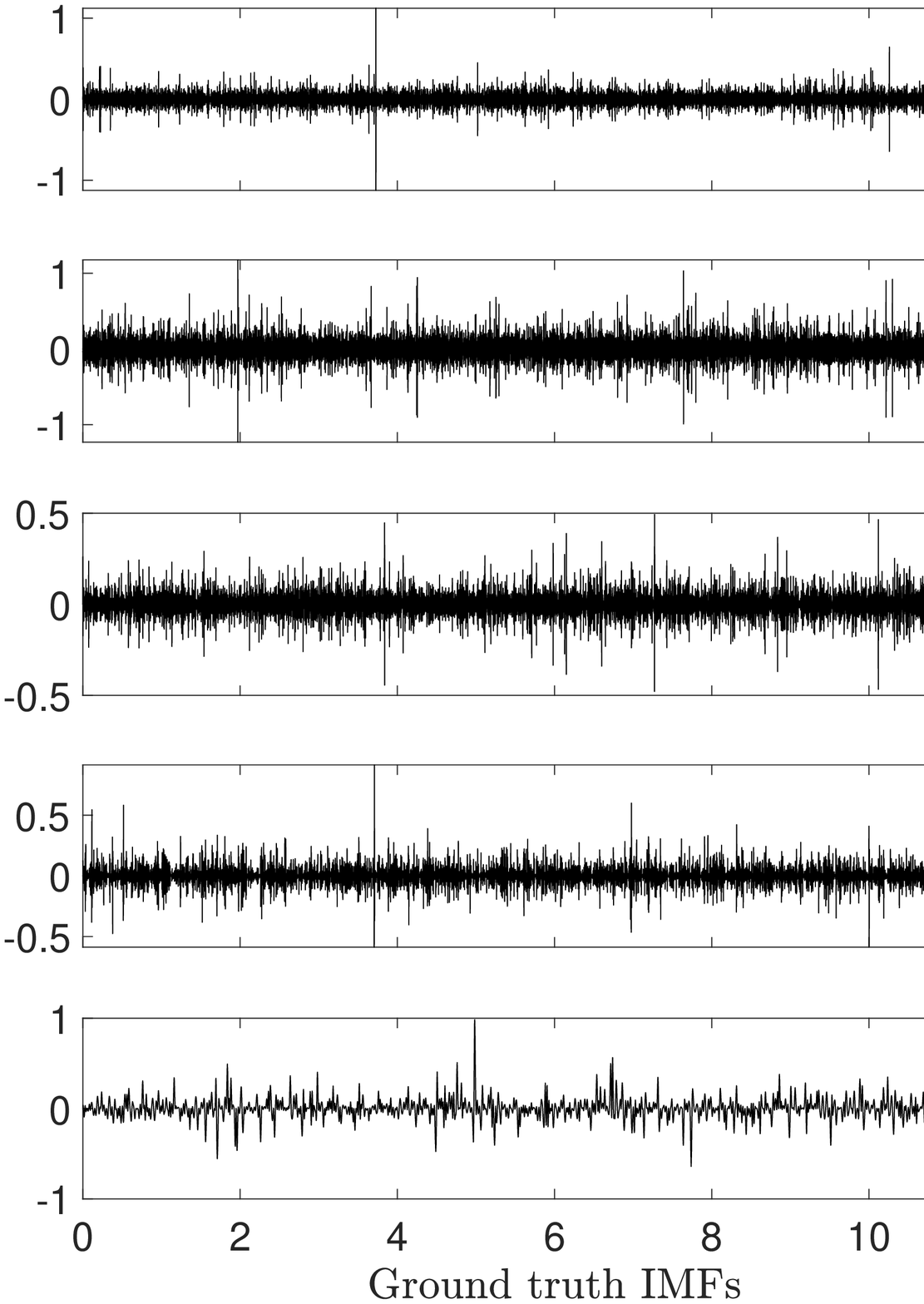} }}%
    ~
    \subfloat{{\includegraphics[width=0.48\textwidth]{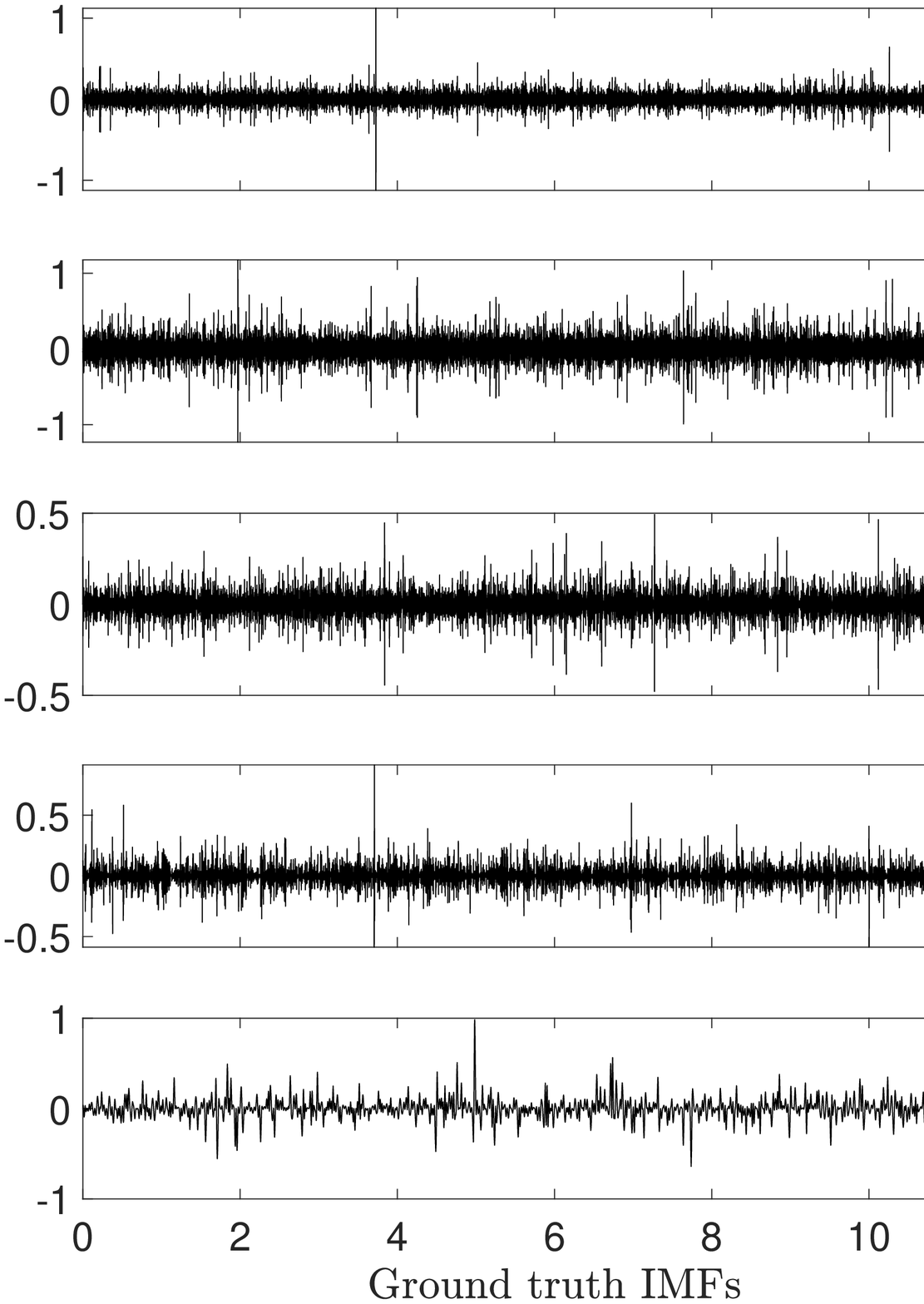} }}%
    \caption{ First and third column from the left, ground truth IMFs. Second and fourth column, the differences between the ground truth signal and the corresponding IMF produced using the FIF and dFIF algorithm, respectively.}\label{fig:Test_ex_1_diff}
\end{figure}

We consider now a toy example containing almost twelve millions sample points, produced using five oscillatory signals which are plotted in the first column of Figure \ref{fig:Test_ex_1_diff}.

\begin{figure}
    \centering
    \subfloat{{\includegraphics[width=0.48\textwidth]{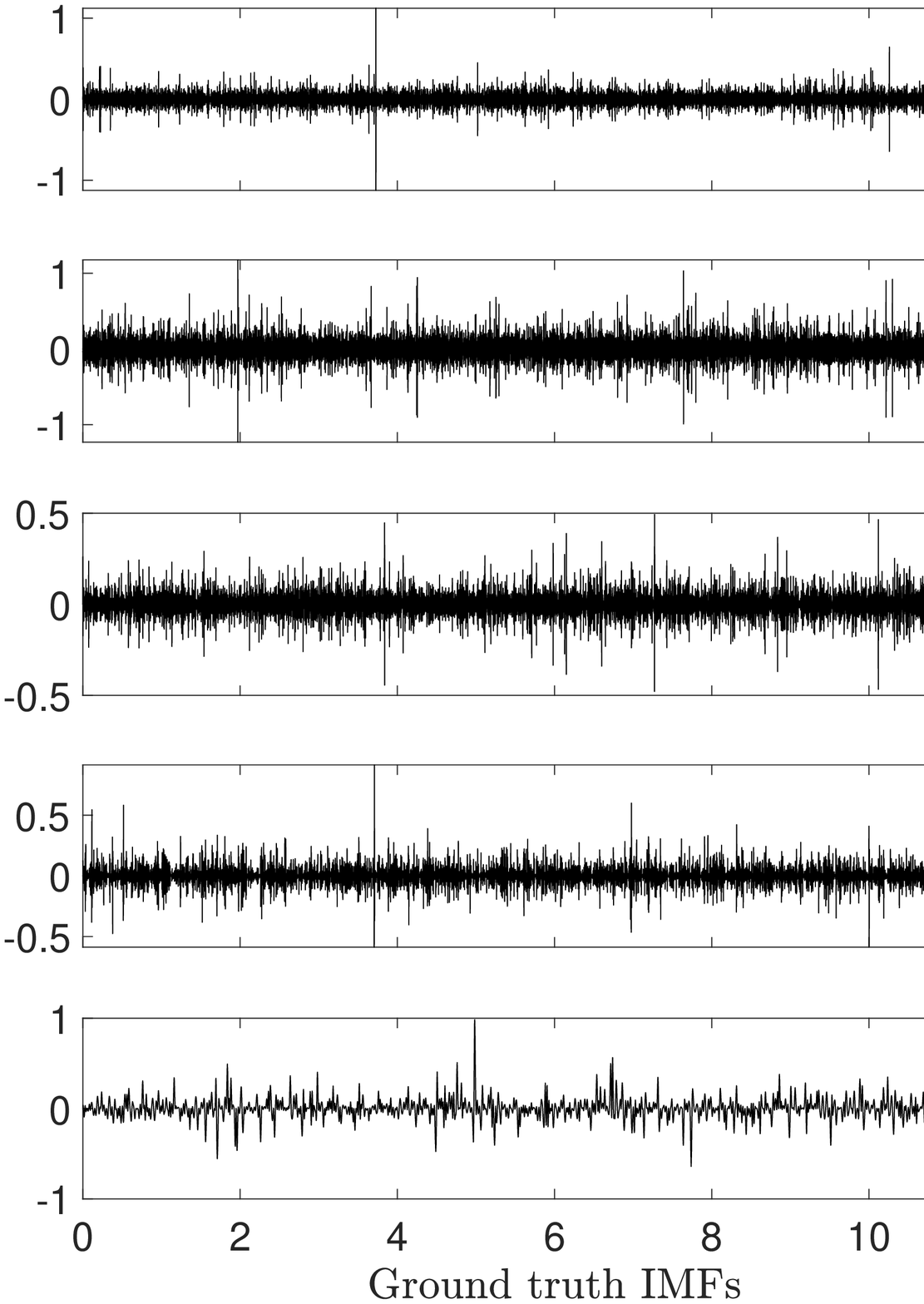} }}%
    ~
    \subfloat{{\includegraphics[width=0.48\textwidth]{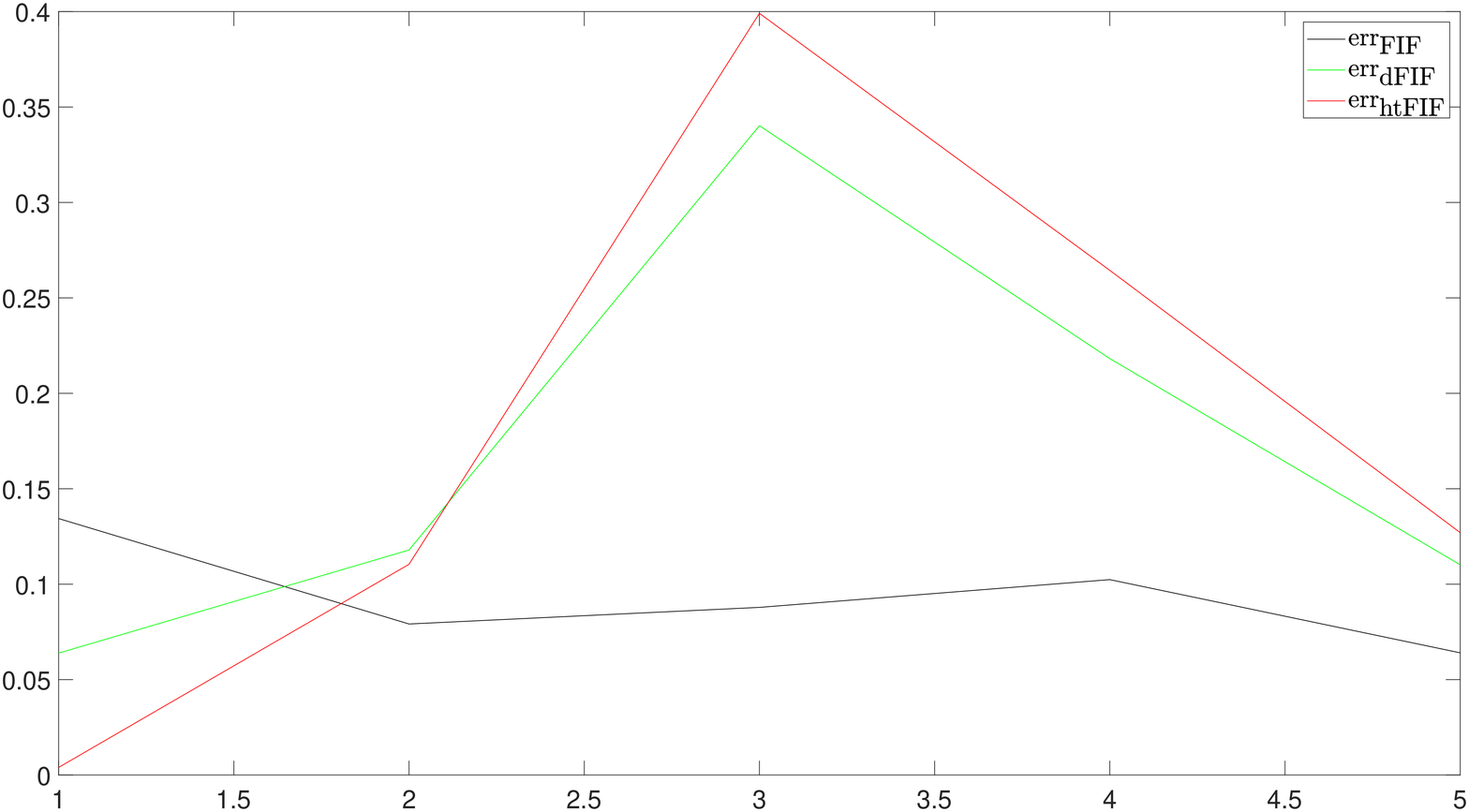} }}%
    \caption{ First column from the left, ground truth IMFs. Second column, the differences between the ground truth signal and the corresponding IMF produced using the htFIF algorithm. Right panel relative errors with respect to the ground truth. In black the relative error for each IMF produced with the FIF algorithm. In green and red the ones relative to the dFIF and htFIF decompositions, respectively.}\label{fig:Test_ex_1_diff2}
\end{figure}

We first run the decomposition using the standard FIF\footnote{We set \textbf{alpha} equal to \textbf{ave}, \textbf{Xi} $=2$, \textbf{delta} $=0.01$.} method. The produced IMFs are compared with the ground truth components using both a point by point subtraction, second column in Figure \ref{fig:Test_ex_1_diff}, and the relative error curve, right panel in Figure \ref{fig:Test_ex_1_diff2}.

Then we apply dFIF and htFIF methods and make similar comparisons with the ground truth data, as shown in the second and fourth column of Figure \ref{fig:Test_ex_1_diff} and in the right panel of Figure \ref{fig:Test_ex_1_diff2}.
As expected, on average, both methods produce less accurate decompositions. Whereas, from a computational time point of view, Table \ref{tab:Ex_2}, we clearly have an improvement with both newly proposed techniques.
\begin{center}
\begin{tabular}{c||c|c|c|c}
  Algorithm & IF & FIF & dFIF & htFIF \\
  \hline
  \hline
  Time (s) & $\quad 3329.7 \quad$ & $\quad 40.72 \quad$ & $\quad 26.04 \quad$ & $\quad 24.23 \quad$ \\
\end{tabular}
\captionof{table}{Computational time for Example 2}\label{tab:Ex_2}
\end{center}

\subsection{Example 3}
\begin{figure}%
    \centering
    \subfloat{{\includegraphics[width=0.48\textwidth]{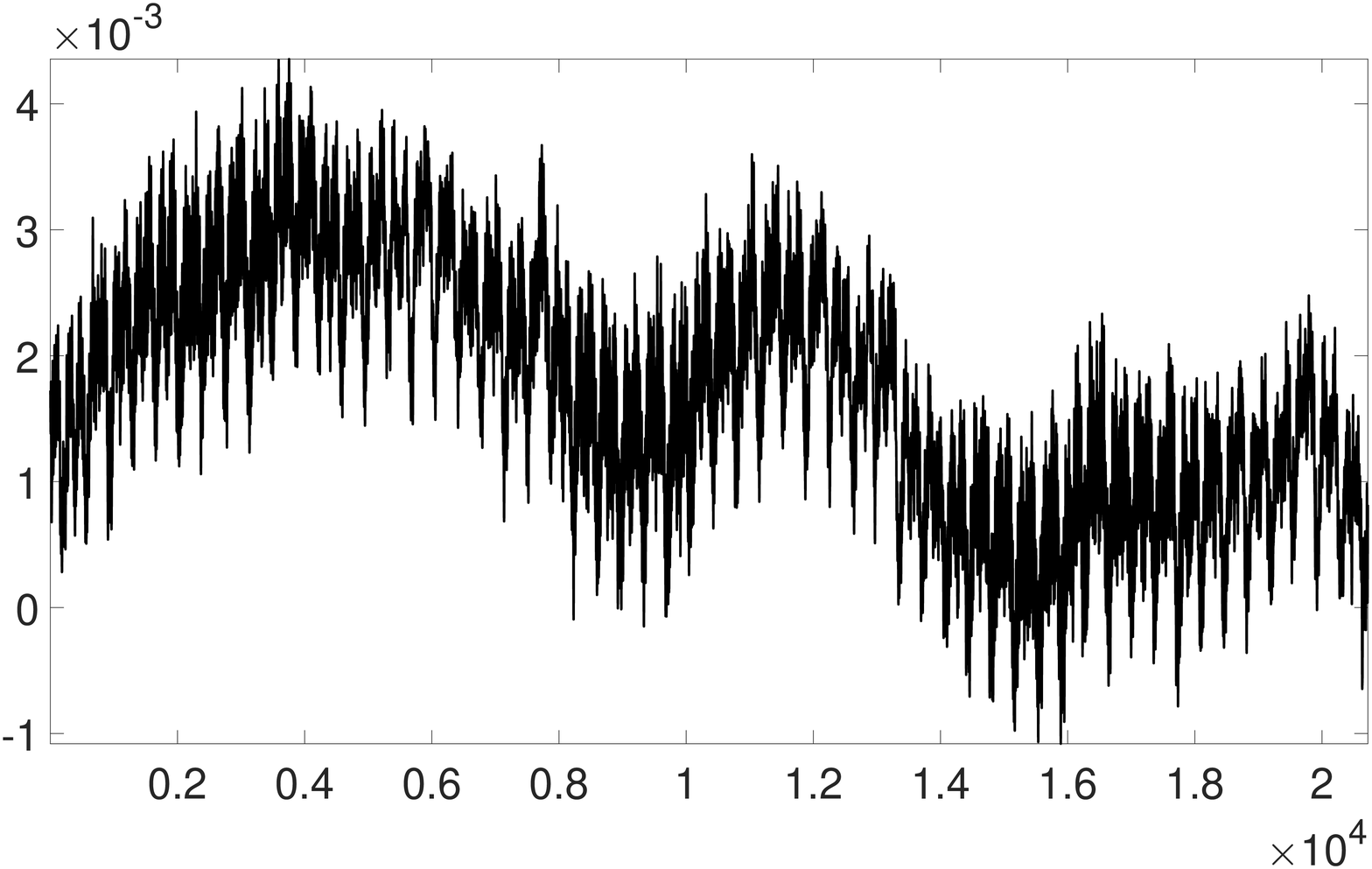} }}%
    ~
    \subfloat{{\includegraphics[width=0.48\textwidth]{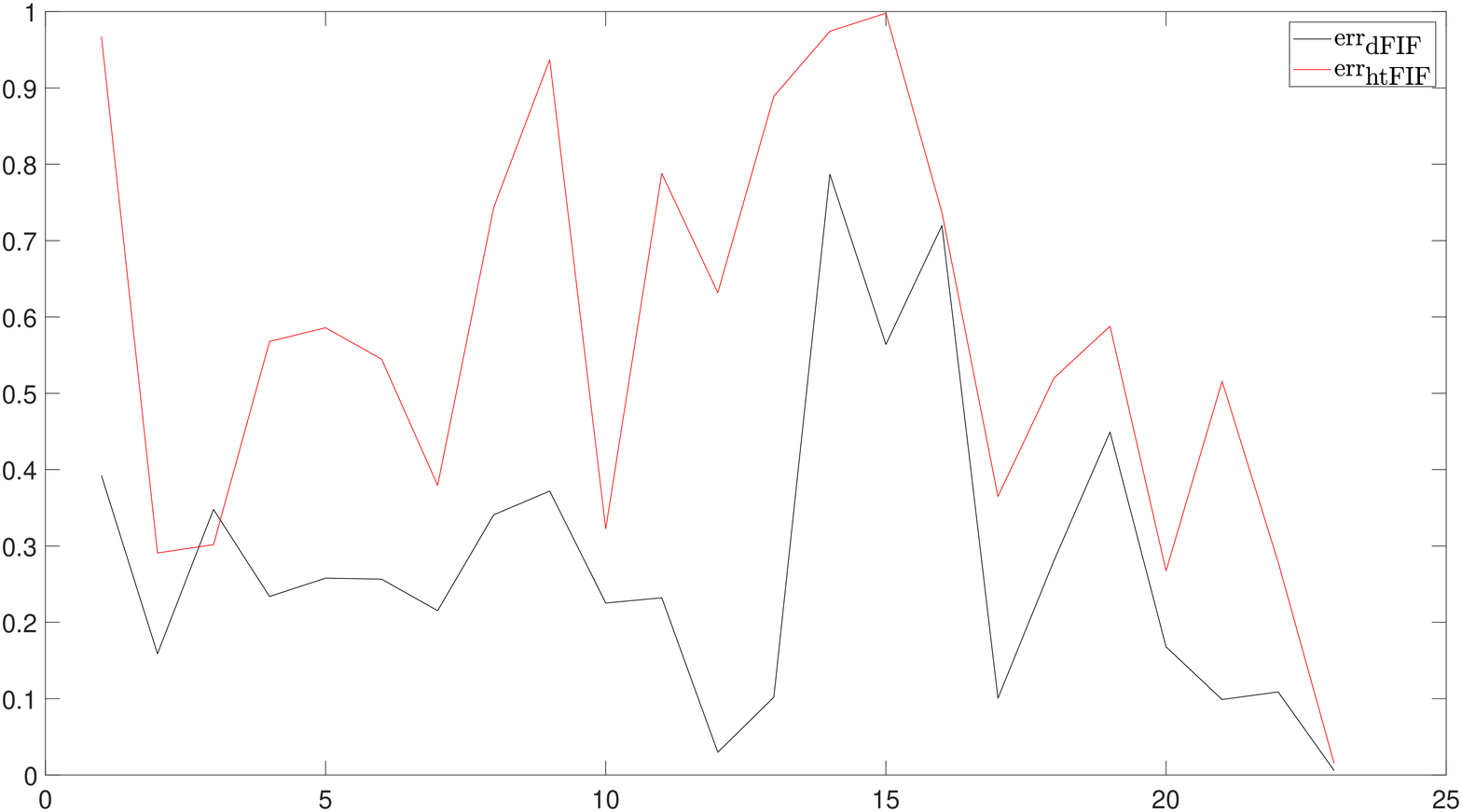} }}%
        \caption{ Left panel, Length of the day signal from 1-1-1962 to 30-9-2018. Right panel, the relative errors $\textrm{err}_{\dFIF}$ and $\textrm{err}_{\htFIF}$ computed for each IMF.}\label{fig:Test_ex_3}
\end{figure}

We test now the three algorithms on a real life signal, the Length Of the Day (LOD) dataset\footnote{LOD dataset is maintained by the The International Earth Rotation and Reference Systems Service and it can be downloaded from \url{http://hpiers.obspm.fr/eoppc/eop/eopc04/eopc04.62-now}. A guide describing how the dataset has been generated can be downloaded from \url{http://hpiers.obspm.fr/eoppc/eop/eopc04/C04.guide.pdf}}, which is shown in the left panel of Figure \ref{fig:Test_ex_3}. The LOD measures the fluctuations of the daily rotation time length of the Earth. This signal has been extensively studied in the past, \cite{Huang2008geophysical} and references therein, and its subtle variations have been related to interactions between the dynamic atmosphere and Earth itself.

\begin{figure}%
    \centering
    \subfloat{{\includegraphics[width=0.48\textwidth]{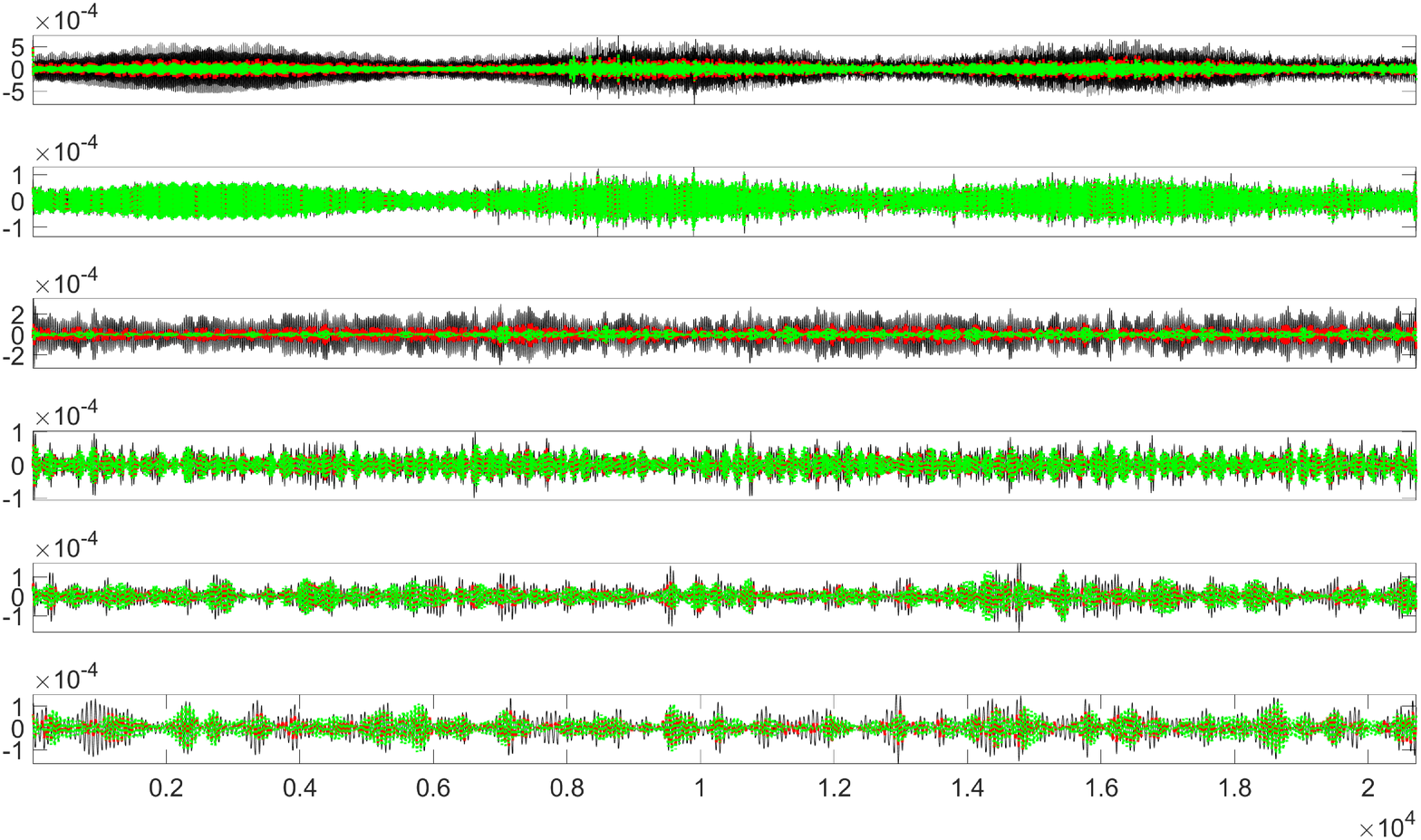} }}%
    ~
    \subfloat{{\includegraphics[width=0.48\textwidth]{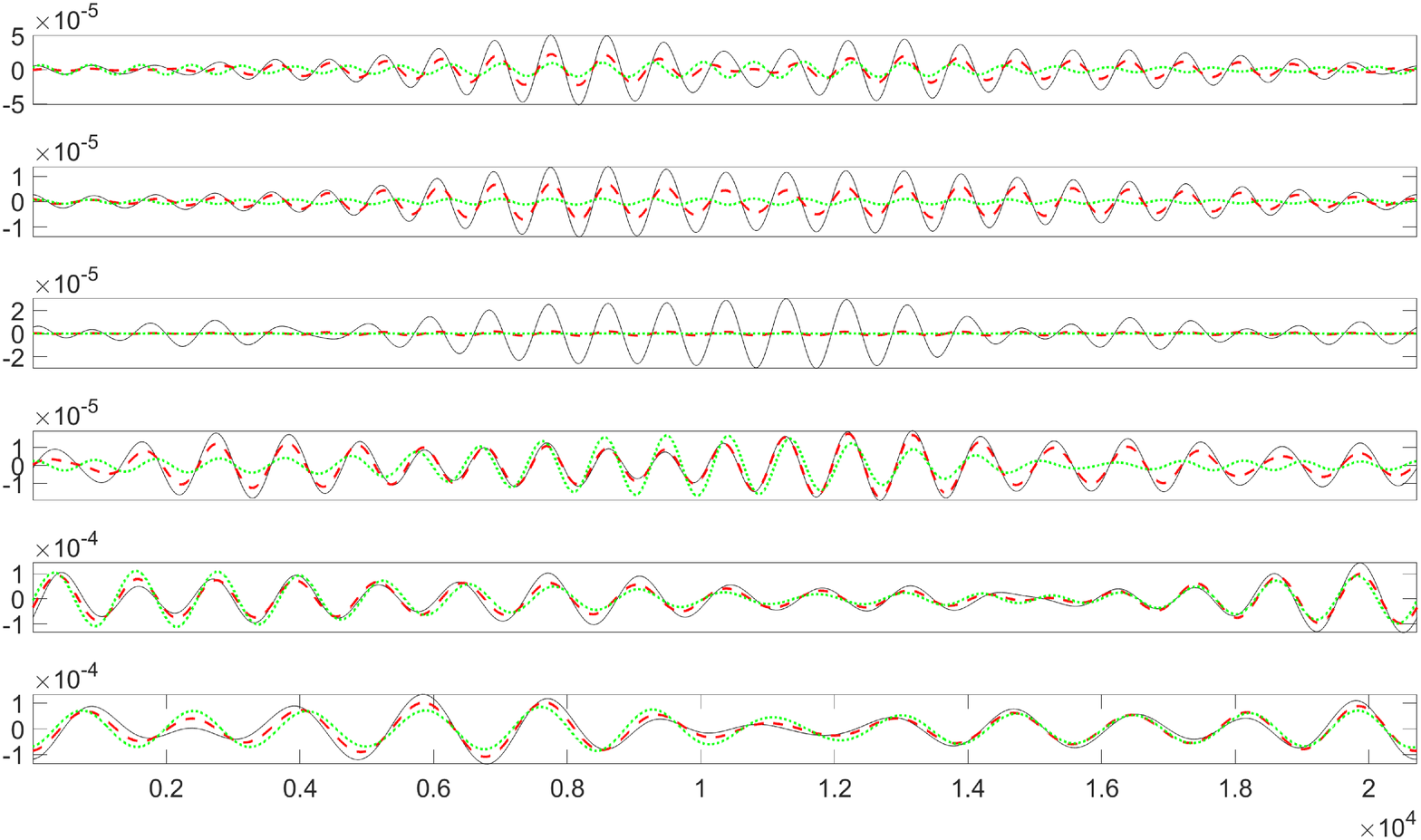} }} \\
    \subfloat{{\includegraphics[width=0.48\textwidth]{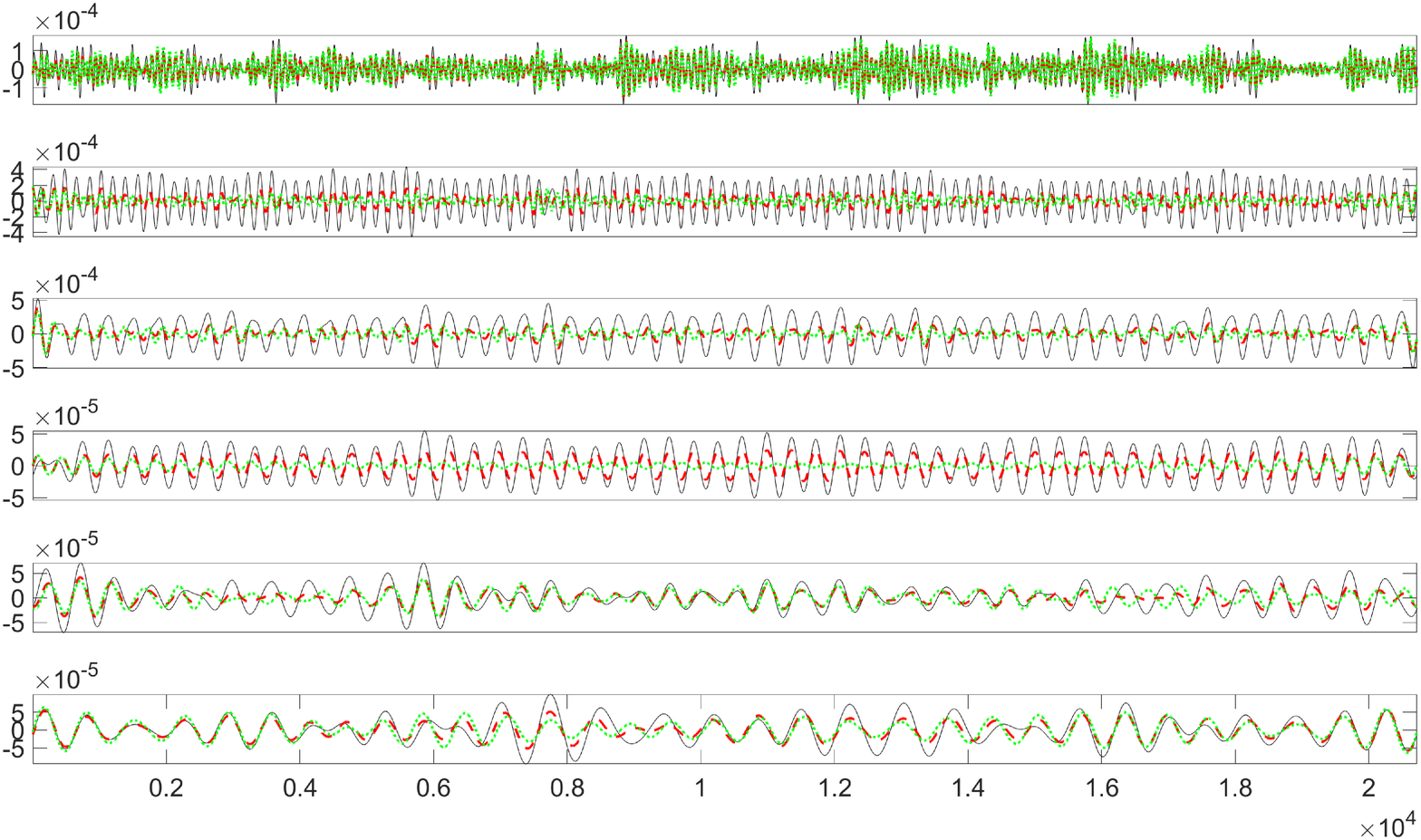} }}%
    ~
    \subfloat{{\includegraphics[width=0.48\textwidth]{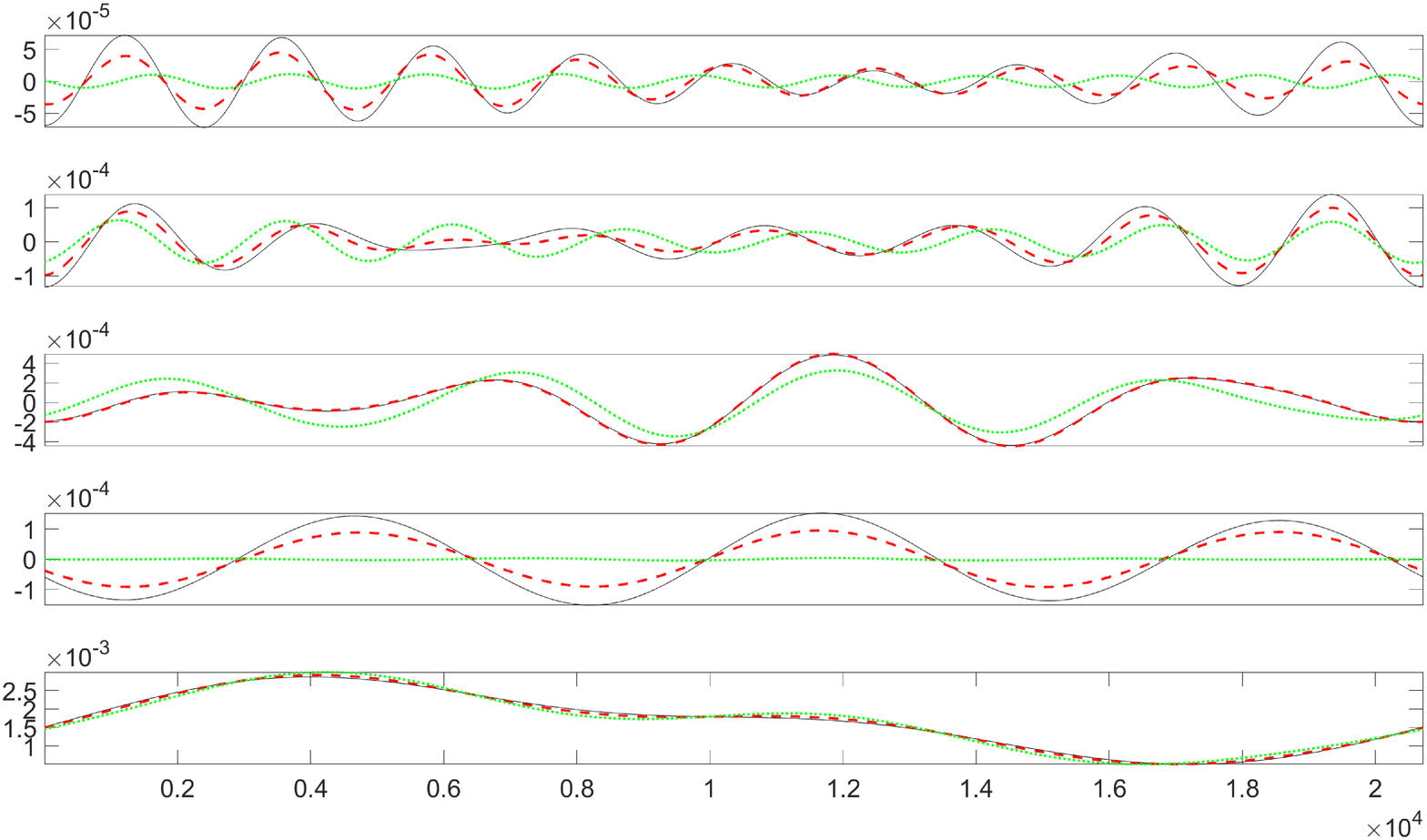} }}%
    \caption{IMFs produced by FIF, dFIF and htFIF methods ordered from highest frequencies, top left corner, to the trend, bottom right. In solid black the FIF IMFs, in dashed red the dFIF IMFs, and in green dots the htFIF IMFs.}\label{fig:Test_ex_3_comparisons}
\end{figure}

In Figure \ref{fig:Test_ex_3_comparisons} we compare the decomposition produced using the FIF\footnote{We set
\textbf{alpha} equal to \textbf{Almost\_min}, \textbf{Xi} $=1.8$, \textbf{delta} $=0.001$.} algorithm with the ones produced using dFIF and htFIF. To be able to compare properly the IMFs produced by the three methods we computed each IMF using the remainder generated by the FIF algorithm. In doing so we can measure meaningful relative errors $\textrm{err}_{\dFIF}$ and $\textrm{err}_{\htFIF}$ for each IMF, right panel of Figure \ref{fig:Test_ex_3}.

The computational time required by the three methods as well as the basic IF are reported in Table \ref{tab:Ex_3}. 
The results regarding the dFIF and htFIF are compatible with the fact that on average in this example the FIF method required 17 iterations to compute each of the 22 IMFs produced in the decomposition, whereas the two newly proposed methods are direct therefore they do not need any iteration.

We mention that, as expected from the theory, the decompositions produced with IF and FIF are identical up to machine precision.

\begin{center}
\begin{tabular}{c||c|c|c|c}\label{tab:Ex_3}
  Algorithm & IF & FIF & dFIF & htFIF \\
  \hline
  \hline
  Time (s) & $\quad 38.09 \quad$   & $\quad 0.47 \quad$ & $\quad 0.069 \quad$ & $\quad 0.057\quad$ \\
\end{tabular}
\captionof{table}{Computational time for Example 3}\label{tab:Ex_3}
\end{center}

From the relative errors plot we can see that in this example both newly proposed methods produce decompositions that are in general different from the FIF one, as we do expect, and that the dFIF method performs always better than the htFIF algorithm.

The relative errors confirm what we expected. The htFIF requires just one parameter tuning, but its decomposition appears to be less similar to the original FIF method. The dFIF method requires two parameters tuning, but it produces a decompositions much closer to the FIF one.
We remind once more that the FIF decomposition, that we are using here as a reference, is not producing in general the exact ground truth as shown in the second example. In fact it requires by itself the tuning of the \verb"delta" parameter.

\section{Conclusions and Outlook}\label{sec:Conclusions}

In this work we propose two alternative reformulations of the so called Iterative Filtering method. These two new versions of the algorithm transform the iterative original method into two direct techniques. They both allow to reduce the computational time required by the original algorithm in the decomposition of signals. The numerical examples proposed in this work confirm this ability and they show that the proposed algorithms can be used as first approximation methods that can quickly decompose long data sets.


\end{document}